\documentclass[12pt]{amsart}
\usepackage{amssymb}

\usepackage{euscript}

\let\cal\mathcal

\makeatletter \@mparswitchfalse \makeatother

\def\defeq{\buildrel{\mathrm{def}}\over=}

\newtheorem{theorem}{Theorem}[section]

\newtheorem{lemma}[theorem]{Lemma}
\newtheorem{corollary}[theorem]{Corollary}
\newtheorem{proposition}[theorem]{Proposition}

\theoremstyle{definition}
\newtheorem{remark}{Remark}

\newtheorem{remarks}{Remarks}

\newtheorem{example}{Example}
\newtheorem{ack}{Acknowledgments}

\newcommand{\bbN}{{\mathbb{N}}}
\newcommand{\bbC}{{\mathbb{C}}}

\newcommand{\al}{\alpha}

\newcommand{\g}{\gamma}

\newcommand{\la}{\lambda}

\newcommand{\s}{\sigma}
\newcommand{\f}{\varphi}

\newcommand{\La}{\Lambda}

\newcommand{\G}{\Gamma}

\newcommand{\Span}{\operatorname{span}}
\newcommand{\supp}{\operatorname{supp}}
\newcommand{\card}{\operatorname{card}}

\newcommand{\sgn}{\operatorname{sgn}}

\newcommand{\disp}{\displaystyle}
\newcommand{\lb}{\label}

\newcommand{\lra}{\longrightarrow}

\newcommand{\wtw}{if and only if }

\def\N{\mathbb N}
\def\C{\mathbb C}

\def\Z{\mathbb Z}

\def\R{\mathbb R}
\let\<\langle

\numberwithin{equation}{section}

\begin{document}

% bring back \eqalign from plain TeX
%\def\eqalign#1{\null\,\vcenter{\openup\jot
% \ialign{\strut\hfil$\displaystyle{##}$&$\displaystyle{{}##}$\hfil
%    \crcr#1\crcr}}\,}

\let\\\cr
\let\phi\varphi
\let\union\bigcup
\let\inter\bigcap
\def\supp{\operatorname{supp}}
\def\vsupp{\operatorname{v-supp}}
\def\sup{\operatorname{sup}}
\def\inf{\operatorname{inf}}
\def\dim{\operatorname{dim}}
\def\Span{\operatorname{span}}
\def\card{\operatorname{card}}
\def\sgn{\operatorname{sgn}}
\def\log{\operatorname{log}}
\def\max{\operatorname{max}}
\def\min{\operatorname{min}}
\let\emptyset\varnothing

%{\today}, Apr. 25

\title{${\bf 1}$-complemented subspaces of  spaces with ${\bf 1}$-unconditional bases}
\author{Beata Randrianantoanina}
\address{Department of Mathematics \\ The University  of Texas at Austin\\
Austin, TX 78712}
\curraddr{Mathematical Sciences Research Institute\\
1000 Centennial Drive\\
Berkeley, CA 94720}
\email{brandri@math.utexas.edu}
\thanks{This research was started while the author participated in 
the  Workshop on Linear Analysis and Probability held at Texas A\&M 
University in College Station, Texas, supported by the NSF Research 
Group grant MCS DMS-9311902, and completed during the stay at Mathematical Sciences Research Institute in Berkeley, California, supported by the 
NSF  grant   DMS-9022140.}

\mbox{       }

%\centerline{Preliminary version}

%page 1

\begin{abstract}
We prove that if $X$ is a complex strictly monotone sequence
space with $1$-unconditional basis, $Y \subseteq X$ has no bands
isometric to $\ell_2^2$ and $Y$ is the range of norm-one projection from
$X$, then $Y$ is a closed linear span a family of mutually
disjoint vectors in $X$.

We completely characterize $1$-complemented subspaces and norm-one
projections in complex spaces $\ell_p(\ell_q)$ for $1 \leq p, q < \infty$.

Finally we give a full description of the subspaces that are spanned
by a family of disjointly supported vectors and which are
$1$-complemented in (real or complex) Orlicz or Lorentz sequence spaces. In particular if an Orlicz or
Lorentz space $X$ is not isomorphic to $\ell_p$ for some $1 \leq p <
\infty$ then the only subspaces
of $X$ which are $1$-complemented and disjointly supported are the 
closed linear spans of  block bases with constant
coefficients.
\end{abstract}
\maketitle

\section{Introduction}

%page 3

Projections and norm one projections have been studied by many authors. The question about the form of a (norm-one) projection and the structure of 
its range  arises naturally not only in geometry
of Banach spaces, but also in approximation theory, spectral theory,
ergodic theory; see, e.g., the surveys \cite{ChP,Doust} for more detailed
discussions of applications.

The difficulty in studying $1$-complemented subspaces  of spaces with
$1$-unconditional bases arises from the following classical fact due to Lindenstrauss \cite{L72} (cf. also \cite[Theorem~3.b.1]{LT1})

\begin{theorem}\lb{nohope}
Every space $Y$ with a
$1$-unconditional basis is $1$-complemented in some symmetric space $X$.
\end{theorem}

Thus it seems hopeless to give any characterization of $1$-complemented 
subspaces of, even symmetric, spaces with
$1$-unconditional bases.

The only class of spaces where the full characterization of $1$-complemented subspaces was available are  the classical spaces $\ell_p$ and $c_0$. 
Namely it is well
known that every subspace of a Hilbert space is $1$-complemented (with
the unique orthogonal projection) and in $\ell_p$, for $p \neq 2, \infty$, 
we have
the following result:
%page 4
\begin{theorem}[{\cite{Ando}, cf.\ also \cite[Theorem 2.a.4]{LT1}}]
\label{lp}
Let $F \subset \ell_p$, where $1 \leq p < \infty$, $p \neq
2$. Then $F$ is $1$-complemented if and only if
\begin{itemize}
\item[(a)] $F$ is isometric to $\ell_p^{\dim F}$, 
\end{itemize}

\noindent
or
\begin{itemize}
\item[(b)] $F$ is spanned by a family of mutually disjoint vectors.
\end{itemize}
\end{theorem}

It is clear that Theorem \ref{lp}(a) cannot be extended to other
spaces. Namely Lindberg \cite{L73} showed a class of Orlicz functions
$\phi$ (for necessary definitions see Section~\ref{pre}) so that
there exists a $1$-complemented subspace $F$ in $\ell_\phi$ such that $F$
is not even isomorphic to $\ell_\phi$. Altshuler, Casazza and B.~L.~Lin
\cite{ACL} showed a similar example in the class of Lorentz sequence
spaces $\ell_{w,p}$. %
%page 5
However, both of these examples were spanned by a family of mutually
disjoint vectors; in fact they were closed linear spans of a block
basis with constant coefficients. 
Also the symmetric space $X$ constructed in Theorem~\ref{nohope} was 
such that $Y$ was isometrically isomorhic to a closed linear span of a block
basis with constant coefficients. 

It is well known that all such spans
are $1$-complemented in any symmetric space ({\cite[Theorem 3.a.4]{LT1}}),
so in fact all of those examples satisfy condition $(b)$ of Theorem~\ref{lp}.

In this paper we prove that indeed Theorem~\ref{lp}(b) can be extended
to a large class of $1$-complemented subspaces of complex spaces with
$1$-unconditional basis. %
%page 6

Namely, if $X$ is a complex, strictly monotone sequence space with a
$1$-unconditional basis, $Y \subset  X$ is $1$-complemented in $X$, and
$Y$ does not contain a band isometric to $\ell_2^2$, then $Y$ is spanned
by a family of disjointly supported vectors (see Corollary~\ref{cor}). 
It is clear that our restrictions on $X$ and $Y$ are necessary (see Remark 
after Corollary~\ref{cor} and examples in Section~\ref{secl2}).

The above-mentioned assumption on $Y$ is satisfied, for example, in all spaces
$X$ that do not have a $1$-complemented subspace isometric to $\ell_2^2$. We
discuss it in greater detail in Section~\ref{secl2}.

In Theorem \ref{main} we also describe the form of general
$1$-complemented subspaces of complex strictly monotone spaces.

Our method of proof %
%page 7
cannot be extended to real sequence spaces. We use in particular the
fact that
every $1$-complemented subspace of a complex space with
$1$-unconditional basis also has a $1$-unconditional basis. The
analogous fact is false in real spaces \cite{Le79,BFL} (see \cite{R86} for the disscusion in special real spaces).

As a consequence of Theorem \ref{main} we obtain a complete
characterization of $1$-comple\-mented subspaces of complex $\ell_p(\ell_q)$, where $1<p, q < \infty$ (Theorems \ref{lplq} and \ref{l2lq}). %
%page 8

Further we study the subspaces that are spanned by disjointly
supported vectors and are $1$-complemented in $X$. Calvert and
Fitzpatrick \cite{CF86} showed that if all disjointly supported
subspaces are $1$-complemented in $X$ then $X$ is isometric to $\ell_p$,
for some $p$, $1 \leq p < \infty$, or to $c_0$.

 In Section~\ref{gen} we
completely characterize the disjointly supported subspaces that are
$1$-complemented in Orlicz and Lorentz sequence spaces
(Theorems \ref{orlicz} and \ref{glorentz}). In particular, if a
Lorentz or Orlicz
space $X$ is not isomorphic to $\ell_p$ for some $1 \leq p < \infty$ then
the only disjointly supported subspaces that are $1$-complemented are
those guaranteed by \cite[Theorem 3.a.4]{LT1}, i.e. spanned by a block
basis with constant coefficients. The results of
Section~\ref{gen} are valid for both real and complex spaces. %
%page 9
%\filbreak

\begin{ack}
I wish to express my gratitude to Professors W. B. Johnson and A. Koldobsky for  valuable suggestions.
\end{ack}

\section{Preliminaries}
\label{pre}

In the following we will consider complex Banach spaces $X$ with a
normalized $1$-unconditional basis $\{e_i\}_{i \in I}$,  where $\N
\supseteq I = \{1, \dots , {\dim X}\}$. Our results are valid in both
the finite- and infinite-dimensional case.

If $x \in X$ we will write $x=(x_i)_{i \in I}$ if
$$
x = \sum^{\dim X}_{i \in 1} x_i e_i \quad \text{ and } \quad \supp x = \{i
\in \N: x_i \neq 0\}.
$$
%page 10
For $x \in X$ we will denote by $x^{*}$ (or sometimes by $x^{N}$)
a \emph{norming functional} for $x$, that is,  $x^{*} \in X^{*}$,
$\|x^{*}\|_{X^{*}} = 1$ and $ x^{*}(x) = \|x\|_X$.

We say (following \cite{KW}, see also \cite{ST73}) that an element $x\in X$ is {\it hermitian} if there exists a hermitian projection $P_x$ from $X$ onto $\Span\{x\}$.

Equivalently, $x$ is hermitian if and only if for all $y \in X$,
$y^{*}$ norming for $y$, and $x^{*}$ norming for $x$ we have
$$
x^{*}(y)y^{*}(x) \in \R.
$$

The set of  all hermitian elements is denoted $h(X)$.
%page 11
%1.5 hour so far

Let $\{H_\la\ : \la\in\La \}$ be the collection of maximal linear subspaces of $h(X)$.
Then $\{H_\la\ : \la\in\La \}$ are called {\it Hilbert components} of $X$. Kalton and
Wood \cite{KW} proved that Hilbert components are well-defined and mutually disjoint.

A Hilbert component $H_\la$ is called {\it nontrivial\/} if $\dim H_\la >1$.

For the careful analysis and properties of Hilbert components of various 
spaces we refer to \cite{KW} and to expository papers \cite{F83,R83}. 
Here we just want to recall some properties which will be used in our 
arguments.

First, recall that if $X$ has 1-unconditional basis $\{e_i\}_{i\in I}$ 
then each
basis element is hermitian. Moreover Kalton and Wood proved the following:

\begin{theorem}[{\cite[Theorem 6.5]{KW}}]
\label{KW}

Let $X$ be a Banach space with a normalized $1$-uncon\-ditional basis.
Then $x \in X$ is hermitian in $X$ if and only if
\begin{itemize}
\item[(i)]
$\|y\|_X = \|y\|_2$ for all $y \in X$ with $\supp y \subset \supp x$,
and
\item[(ii)] for all $y, z \in X$ with $\supp y \cup \supp z \subset
\supp x$ and for all $v \in X$ with $\supp v \cap \supp x =
\varnothing$ if $\|y\|_X = \|z\|_X$ then $\|y+v\|_X = \|z+v\|_X$.
\end{itemize}
\end{theorem}
%page 12

For our main result we will need the following two facts.

\begin{proposition}[{\cite[Lemma~5.2]{KW}}] \lb{kw}
Suppose that $x, y$ are hermitian elements in $X$. Denote by $x^*$ a norming functional
for $x$.

If $x^*(y)\neq 0$ then $\Span \{x, y\} \subset h(X)$.
\end{proposition}

\begin{proposition}[{\cite[Lemma~4]{F83}}] \lb{yherm}
Suppose that $X$ has 1-unconditional basis $\{e_i\}_{i\in I}$ and let 
$P: X\lra X$
be norm one projection with range of $P$ equal to $Y$. Then for all $i\in I$, $Pe_i$
is a hermitian element  in $Y$.
\end{proposition}

We will also frequently use the following well-known fact:
\begin{proposition}
\label{conj}
Let $X$ be a Banach space with a $1$-unconditional basis. Suppose that
$Y \subset X$ is $1$-complemented and a norm-one projection $P:X
\to Y$ is given by
$$
P(x) = \sum_i y_i^{*}(x) y_i
$$
where $Y = \overline\Span\{y_i\}$ and $y_i^{*}$ is norming for $y_i$ for
all $i$.

Then for any $y \in Y$ there exists $y^{*}$ norming for $y$ and
constants $K_i$ so that
$$
y^{*} = \sum_i K_i y_i^*.
$$
\end{proposition}

Moreover, we have
\begin{proposition}[Calvert \cite{Cal75}]
\label{calvert}
Let $X$ be a strictly convex reflexive Banach space with strictly convex dual $X^*$. Let $J:X\lra X^*$ be the duality map; 
$\|Jx\| = \|x\|, \ Jx(x) = \|x\|^2.$

Then a closed linear subspace $Y$ of $X$  is $1$-complemented in $X$ \wtw
$J(Y)$ is a linear subspace  of $X^*$.
\end{proposition}

%page 13

Finally we recall a few definitions (see \cite{LT1}).

We say that a Banach space $X$ with $1$-unconditional basis is \emph{
strictly monotone\/} if $\|x+y\| > \|x\|$ for all $x, y \geq 0$ with
$y \neq 0$.

An \emph{Orlicz function} $\phi$ is a convex non-decreasing function
$\phi: [0,\infty)\longrightarrow [0,\infty]$ with $\phi(0) = 0$ and
$\phi(1) = 1$ or $\infty$. To any Orlicz function $\phi$ we associate
the \emph{Orlicz space} $\ell_\phi$ of all sequences of scalars $x =
(x_i)_i$ such that
$$
\sum^\infty_{i=1} \phi \Bigl(\frac{|x_i|}{\rho}\Bigr) < \infty \quad
\text{ for some } \rho > 0,
$$
with the norm
$$
\|x\|_\phi = \inf\biggl\{\rho > 0 : \sum^\infty_{i=1} \phi
\Bigl(\frac{|x_i|}{\rho}\Bigr) < 1\biggr\}.
$$

Let $1 \leq p < \infty$ and let $w=\{w_i\}_{i \in I}$, where $I=\N$
or $I = \{1,\dots, d\}$, be a %
%page 14
non-increasing sequence such that $w_1=1$ and $w_i \geq  0$ for all
$i$. The Banach space of all sequences of scalars $x = (x_i)_{i \in
I}$ for which
$${\disp 
\|x\|_{w,p} = \sup\limits_{\sigma \in \cal{P}(I)}\biggl(\sum_{i \in I}
|x_{\sigma(i)}|^p w_i\biggr)^{\frac{1}{p}} < \infty,}
$$
where $\cal{P}(I)$ is the set of all permutations of $I$, is
denoted $\ell_{w,p}$ and it is called a \emph{Lorentz sequence space}
(another notation frequently used in the literature is $d(w,p))$.

Notice that
\begin{itemize}
\item $\ell_\phi$ is strictly monotone if and only if $\phi(t) > 0$ for
all $t > 0$ and $\phi(t) < \infty $ for all $t \leq 1$.

\item $\ell_{w,p}$ is strictly monotone if and only if $w_i > 0$ for
all $i \in I$.

\item $\ell_p(\ell_q)$ is strictly monotone if and only if $p,q \neq
\infty$.
\end{itemize}
%page 15

For any $ 1 \leq p, q \leq \infty$ we denote be $\ell_p(\ell_q)$ the space
of sequences of scalars $x = (x_{ij})_{i \in I, j \in J}$ such that
$$
\|x\|_{\ell_p(\ell_q)} = \bigl\| ( \| (x_{ij})_{j \in J}\|_{\ell_q})_{i \in
I}\bigr\|_{\ell_p} < \infty.
$$

We follow standard notations as defined in \cite{LT1} and this is also
where we refer the reader for all undefined terms. %
%page 16

\section{General form of contractive projections}

We are now ready to present our main theorem.

\begin{theorem}\lb{main}
Suppose that $X$ is a complex  strictly monotone sequence
space with 1-unconditional basis $\{e_i\}$ and $X\ne \ell_2$ and let $P$ be the
projection of norm 1 in $X$.
Let $\{H_\g\ : \g\in\G\}$ be the collection of Hilbert components of $Y = PX$.
Then $H_\g's$ are disjointly supported as elements of $X$.
\end{theorem}

\begin{proof}   %{Proof of Theorem~\ref{main}}
By Proposition~\ref{yherm} all $\{Pe_i\}_{i\in I}$ are hermitian 
elements of $Y$.
Let $i,j$ be such that $Pe_i$ and $Pe_j$ are not in the same Hilbert 
component of $Y$.
Assume that there exists $k$ such that $k\in \supp Pe_i\cap \supp Pe_j$.
If $Pe_k\neq 0$, then $P^*e_k^*$ is a
norming functional for $Pe_k$ in $Y$
and $(P^*e_k^*,Pe_i)=(e_k^*,Pe_i)\neq 0. $ Thus, by Proposition~\ref{kw},  
$Pe_k, Pe_i$ are in the same Hilbert component of $Y$.
Similarly $Pe_k, Pe_j$ are in the same Hilbert component of $Y$. But 
then $Pe_i,
Pe_j$ are in the same Hilbert component of $Y$ contrary to our assumption.

Thus $Pe_k = 0$.

Now suppose that $Pe_i = \sum_{l\in S} \alpha_l e_l +\alpha_k e_k$ for some $\al_l, \al_k \neq 0$,
where $S= \supp Pe_i\setminus \{k\}$.

Since $P$ is a projection and $P(e_k)=0$ we get
\begin{align*}
\sum_{l\in S} \alpha_l e_l +\alpha_k e_k &=Pe_i = P(Pe_i)  = P\biggl( \sum_{l\in S} \alpha_l e_l\biggr)
+ \alpha_k P(e_k) =\\
&=P\biggl( \sum_{l\in S} \alpha_l e_l\biggr) \ .
\end{align*}
Hence $S\ne \emptyset$ and by strict monotonicity of $X$
$$\Big\| P\biggl( \sum_{l\in S} \alpha_l e_l\biggr)\Big\|
= \Big\| \sum_{l\in S} \alpha_l e_l + \alpha_k e_k\Big\|
> \Big\| \sum_{l\in S} \alpha_l e_l \Big\|$$
which contradicts the fact that $\|P\|=1$.

Thus if $Pe_i,Pe_j$ are not in the same Hilbert component of $Y$ then they are disjoint.
\end{proof}

\begin{corollary}\lb{cor}
Suppose that $X$ is a complex  strictly monotone sequence
space with 1-unconditional basis $\{e_i\}$ and $X\ne \ell_2$ and let $P$ be the projection of norm 1 in $X$.

Suppose that  $Y = PX\subset X$  has no nontrivial Hilbert components. Then
there exist disjointly supported elements $\{y_j\}_{j=1}^m$ ($m= \dim PX \le\infty$)
which span $Y = PX$. Moreover, for all $x\in X$,
$$Px = \sum_{j=1}^m y_j^*(x)y_j,$$
where $\{y_j^*\}_{j=1}^m\subset X^*$ satisfy $\|y_j\| = \|y_j^*\| = y_j^*(y_j) = 1$ for all $j$.
\end{corollary}

\begin{proof}
By Proposition~\ref{yherm} all $\{Pe_i\}_{i\in I}$ are hermitian elements in $Y$.

By our assumption all Hilbert components of $Y$ are one-dimensional so if  $Pe_i,Pe_j$
are linearly independent then they belong to different Hilbert components of $Y$.
Therefore, by Theorem~\ref{main},   if $Pe_i,Pe_j$ are linearly independent then they are disjoint
and
$Y$ can be presented as $\hbox{span}\{Pe_i :i\in I\}$ where $I$ is a
collection of such indices $i,j$ that $Pe_i,Pe_j$ are mutually disjoint.

Then $y_i = {Pe_i}/{\|Pe_i\|}$ for all $i\in I$, and for each $x\in X$ we have
$Px = \sum_{i\in I} C_i y_i$, where $C_i\in \bbC$ are uniquely determined by $x$.
Clearly $y_i^*(x) \defeq  C_i(x)$ satisfies the conclusion of the theorem.

Notice also that $\supp y_i^* = \supp y_i$ for all $i\in I$.
\end{proof}

\begin{remarks}
\begin{enumerate}
\item Notice that the assumption about $X$ being strictly monotone is
important. Indeed, Blatter and Cheney \cite{BCh} (see also \cite{B88}) showed examples of
$1$-complemented hyperplanes in $\ell^3_\infty$ that are not spanned by
disjointly supported vectors.
\item Also the assumption about $Y$ cannot be removed. We
discuss it in greater detail in the next section.
\item As mentioned in Introduction, Calvert and
Fitzpatrick \cite{CF86} showed that if every subspace of the form described in Corollary~\ref{cor} is  $1$-complemented in $X$ then $X$ is isometric to $\ell_p$,
for some $p$, $1 \leq p < \infty$, or to $c_0$.
\end{enumerate} 
\end{remarks}

As a consequence of Corollary~\ref{cor} we can express $1$-complemented subspaces as an intersection of hyperplanes of special simple form (see \cite{BP88} for analogous representation in $\ell_p$).

\begin{corollary}\lb{ker}
Let $X$ and $Y$ be as described in Corollary~\ref{cor}.
Then $Y$ can be presented as intersection of kernels of functionals
$f_j$, such that $\card(\supp f_j)\le 2$ for all $j$.
\end{corollary}

%page 17

\section{$1$-complemented copies of $\ell^2_2$}
\label{secl2}

In this section we discuss in what situation it is possible that a
space $X$ has a $1$-complemented subspace $Y$ with nontrivial Hilbert
components. This clearly reduces to the question of characterizing under
what conditions $X$ can have a $1$-complemented subspace $F$ that is
isometric to $\ell^2_2$. %
%page 18
The question that arises here is:\par\nobreak\medskip
\begin{quote}
\it
Is it possible that a space $X$ with only $1$-dimensional Hilbert
components has a $1$-complemented subspace isometric to $\ell^2_2$?
\end{quote}
\medbreak
One quickly realizes that the answer is yes: 

\begin{example}\lb{exl2lq}
Consider the space
$X= \ell_2(\ell_1)$ and for $x=(x_{ij})_{i,j} \in \ell_2(\ell_1)$ let $Px = (x_{i1})_i$.   Then $PX$ is isometric to
$\ell^2_2$.

Further there exists an orthogonal projection $Q$ of $PX$ onto a span of any collection of orthogonal vectors from $PX = \ell_2$. 
\end{example}

In Lemma \ref{l2cor} we show that if $\ell_p(\ell_q)$
has a $1$-complemented subspace $F$ that is isometric to $\ell^2_2$
then  either $q=2$ and $F$ is contained in a Hilbert component of $\ell_p(\ell_q)$  or $p=2$ and $F$ is similar to the range of $QP$ in Example~\ref{exl2lq}.

 Below we present two more examples: a $1$-complemented
copy of $\ell^2_2$ in a Lorentz space and in a Orlicz space. We do not know a
full characterization of spaces $X$ that have $1$-complemented copies
of $\ell^2_2$, but we suspect that if $X$ is a Lorentz or Orlicz space
then $X$ has to be very similar to the examples presented below. %
%page 19

\begin{example}
Consider Lorentz space $\ell_{w,2}$ with weight $ w = (1,1, w_3, w_4,
\dots)$. Then $\Span \{e_1, e_2\} \subset \ell_{w,2}$ is isometric to
$\ell^2_2$ and clearly it is $1$-complemented.
\end{example}

\begin{example}\lb{ex2}
Consider 4 dimensional Orlicz space $\ell_\varphi$ where
$$\varphi (t) = \left\{
\begin{array}{ll}
t^2&\mbox{if }\ 0\le t\le a\\
(1+a)t-a&\mbox{if }\ a\le t\le1\end{array}
\right.$$
and $\sqrt{2/3} <a<1$.
That is,
$$\|(x_1,\ldots,x_4)\|_\varphi = \inf
\left\{ \lambda : \sum_{i=1}^4 \varphi \left( \frac{|x_i|}{\lambda}\right)
\le 1\right\}\ .$$
In fact we have $\sum_{i=1}^4 \varphi \left( \frac{|x_i|}{\|x\|_\varphi}
\right)=1$ for all $x\in \ell_\varphi$.
Let $F= \ker (e_1^* + e_2^* + e_3^*) = \hbox{span}\{e_1-e_2, e_1-e_3,e_4\}$.
then, if $x= (x_1,x_2,x_3,x_4)\in F$ then
$\|x\|_\varphi$ is a number such that
\begin{equation}\lb{ex21}
\sum_{i=1}^3 \frac{|x_i|^2}{\|x\|_\varphi^2} + \varphi \left( \frac{|x_4|}{
\|x\|_\varphi}\right) = 1
\end{equation}

Indeed,  suppose that $x\in F$ and denote $\|x\|_\varphi =\alpha$.
Then
\begin{equation}\lb{ex22}
1= \sum_{i=1}^4 \varphi\left( \frac{|x_i|}{\alpha}\right)
\ge \sum_{i=1}^4 \frac{|x_i|^2}{\alpha^2}
\end{equation}
Assume that there is $j$, $1\le j\le3$ such that
$\frac{|x_j|}{\alpha} >a$.
Then, by \eqref{ex22}
$$\sum_{i\ne j} \frac{|x_i|^2}{\alpha^2} < 1-a^2\ .$$
Hence
$$\left| \frac{x_1+x_2+x_3}{\alpha}\right|
\ge \frac{|x_j|}{\alpha} - \sum_{\textstyle
\genfrac{}{}{0pt}{}{i=1}{i\ne j}}^3
\frac{|x_i|}{\alpha} \ge a-\sqrt2\ \sqrt{1-a^2} > 0$$
since $a> \sqrt{2/3}$.
But this contradicts the fact that $x\in F = \ker (e_1^* + e_2^* +e_3^*)$.
Thus $\frac{|x_j|}{\alpha} \le a$ for $j=1,2,3$, so
$$1= \sum_{i=1}^4 \varphi \left( \frac{|x_i|}{\alpha}\right)
= \sum_{i=1}^3 \frac{|x_i|^2}{\alpha^2} +\varphi \left(\frac{|x_4|}{\alpha}\right)
\ .$$

Equation \eqref{ex21} and Theorem~\ref{KW} immediately imply that
$e_1-e_2$ and $e_1-e_3$ are hermitian in $F$ and belong to the same Hilbert component
of $F$, but clearly $F$ is not isometric to $\ell_2$.

Moreover \eqref{ex21} implies that $F$ is 1-complemented in $\ell_\varphi$.
Indeed, define $P:\ell_\varphi\to F$ by:
\begin{equation*}
\begin{split}
P(x_1,x_2,&x_3,x_4) = \\ 
&\left( x_1 -\frac13 (x_1+x_2+x_3),
x_2-\frac13 (x_1+x_2+x_3) , x_3 -\frac13 (x_1+x_2+x_3),x_4\right)  .
\end{split}
\end{equation*}
Notice that $Q:\ell_2^3 \to \ell_2^3$ defined by
$$Q(x_1,x_2,x_3) = \left( x_1 -\frac13 (x_1+x_2+x_3),
x_2-\frac13 (x_1+x_2+x_3) , x_3 -\frac13 (x_1+x_2+x_3)\right)$$
is the norm one projection on this Hilbert space.
Let $x\in\ell_\f$. Denote $y=Px$ and $\|Px\|_\varphi =\beta$.
Then, by \eqref{ex21}
\begin{eqnarray*}
1&=&\sum_{i=1}^3  \frac{|y_i|^2 }{\beta^2} + \varphi \left(\frac{|y_4|}{\beta}
	\right) = \frac{\|(y_1,y_2,y_3)\|_2^2 }{ \beta^2}
	+\varphi \left( \frac{|y_4|}{\beta}\right) \\
&=& \frac{\|Q(x_1,x_2,x_3)\|_2^2 }{\beta^2}
	+\varphi\left( \frac{|y_4|}{\beta}\right)
	\le \frac{\|(x_1,x_2,x_3)\|_2^2}{ \beta^2}
	+\varphi \left( \frac{|x_4|}{\beta}\right) \\
&=& \frac{|x_1|^2}{\beta^2} + \frac{|x_2|^2}{\beta^2}
	+ \frac{|x_3|^2}{\beta^2} +\varphi\left( \frac{|x_4|}{\beta}\right)\le\\
&\le& \sum_{i=1}^4 \varphi \left( \frac{|x_i|}{\beta}\right)\ .
\end{eqnarray*}
Thus $\|x\|_\varphi \ge \beta = \|Px\|_\varphi$, i.e. $\|P\|\le1$.
\end{example}

%1.5 hours more
%page 20

We finish this section with a lemma characterizing $1$-complemented
subspaces $F$ in $X$ that are isometric to $\ell^2_2$ in terms of
norming functionals. %
%page 21
\begin{lemma}
\label{l2}

If $\Span \{x,y\} \subset X$ is $1$-complemented in $X$ and
$\Span\{x,y\}$ is isometric to $\ell^2_2$, then there exist
$x^*$ norming for $x$ and $y^*$ norming for $y$ such that,
for all $a, b \in \C$ with $|a|^2 + |b|^2 = 1$, the functional $\bar{a} x^* +
\bar{b} y^*$ is norming for $ax + by$.
\end{lemma}

\begin{proof}%[Proof of Lemma \ref{l2}]
Let $P: X \to \Span\{x,y\}$ be a norm-one projection. Then there exist
$x^*$ norming for $x$ and $y^* $ norming for $y$ such that
$$
P(u) = x^*(u)x + y^*(u)y \quad \text{ for all } u \in X
$$
and $x^*(y) = y^*(x)= 0$.  Moreover, $\Span \{x^*, y^*\} \subset X^*$
is isometric to $\ell^2_2$.  Thus
$$
(\bar{a} x^* + \bar{b} y^*)(ax + by) = a \bar{a} x^*(x) + \bar{b} b
y^*(y)= |a|^2 + |b|^2 = \|\bar{a} x^* + \bar{b} y^*\|_{x^*} \|ax + by\|_x.
\eqno\qed
$$
\let\qed\relax
\end{proof}
%page 22

\section{$1$-complemented subspaces of $\ell_p(\ell_q)$}
\label{seclplq}

The main results of this section (Theorems \ref{lplq} and \ref{l2lq}) completely
characterize $1$-com\-ple\-men\-ted subspaces of the complex space $\ell_p(\ell_q)$, 
$ 1 < p, q < \infty$, $p \neq q$.

To formulate the theorem conveniently we will need some notation.
For all $x$ in $\ell_p(\ell_q)$ we will write
$x=\sum^n_{i=1}x_i$,
where $n = \dim \ell_p \leq \infty$, and $x_i \in \ell_q$ for all
$i \leq n$. We will also write
$$
x= \sum^n_{i=1} \sum^m_{j=1} x_{ij} e_{ij},
$$
where $m = \dim \ell_q \leq \infty$ and $e_{ij}$ are standard
basis elements in $\ell_p(\ell_q)$. %
%page 23
We will distinguish two types of support of $x$:
\begin{itemize}
\item the usual support
$$
\supp x = \{ (i,j) \subset \{1, \dots,n\} \times  \{1, \dots, m\}:
x_{ij} \neq 0\},
$$
which we will sometimes call  the \emph{scalar support} of $x$, and
\item the \emph{vector support}
$$
\vsupp x = \{i \subset \{1, \dots, n\} : \ell_q \ni x_i \neq \vec{0} \}.
$$
\end{itemize}
Thus we will have notions of disjointness of elements in $\ell_p(\ell_q)$ in
the scalar and vector senses. %
%page 24
\begin{theorem}
\label{lplq}
Let $1 < p, q < \infty$ with $p \neq q, 2$. Consider the complex space $\ell_p(\ell_q)$.
 Then $Y \subset \ell_p(\ell_q)$ is
$1$-complemented if and only if there exist $\{v^i\}^{\dim Y}_{i=1} \subset \ell_p(\ell_q)$ so that
$$
Y = \overline\Span \{v^i\}^{\dim Y}_{i=1},
$$
and for all $i\ne j \le \dim Y$ one of the following conditions holds:
\begin{itemize}
\item[(a)] 
 $\vsupp v^i \cap
\vsupp v^j = \varnothing $, 
\end{itemize}

\noindent
or  
\begin{itemize}
\item[(b)] 
 $ \vsupp v^i  =  \vsupp v^j $, $\|v^i_k\|_q=\|v^j_k\|_q$ for all $k\in \vsupp v^i,$ and 
\begin{itemize}
\item[(b1)]
if $q\ne2$ then       $\supp v^i \cap
\supp v^j = \varnothing $, 
\item[(b2)]
 if $q=2$ then $ v^i_k, v^j_k$ are orthogonal for each $k\in \vsupp v^i.$
\end{itemize}
\end{itemize}
\end{theorem}

The structure of 1-complemented subspaces of $\ell_2(\ell_q)$ can be somewhat more complicated as illustrated in Example~\ref{exl2lq}. We have the following:

\begin{theorem}
\label{l2lq}
Let $1 < q < \infty$ with $q \neq  2$. Consider the complex space $\ell_2(\ell_q) = \ell^n_2(\ell^m_q)$, $n,m \in \bbN\cup\{\infty\}$.
 Then $Y \subset \ell_2(\ell_q)$ is
$1$-complemented if and only if there exist $\{v^i\}^{\dim Y}_{i=1} \subset \ell_p(\ell_q)$ so that
$$
Y = \overline\Span \{v^i\}^{\dim Y}_{i=1},
$$
and for all $i\ne j \le \dim Y$ one of the  conditions $(a)$ or $(b)$ of Theorem~\ref{lplq} holds or
\begin{itemize}
\item[(c)] for each $k\in \vsupp v^i \cap
\vsupp v^j$ there exists a constant $C_k\in\bbC$ with
$$v^i_k = C_k v^j_k$$
and for every map $\s:\{1,\dots,n\} \lra \{1,\dots,m\}$ such that $(k,\s(k))\in\supp v^i \cup \supp v^j$ whenever $k\in \vsupp v^i \cup \vsupp v^j$,  the vectors
$$s_\s(v^i) =\left(\|v^i_k\|_q \frac{v^i_{k\s(k)}}{|v^i_{k\s(k)}|}\right)_{k\le n} \in \ell_2,   \text{ and } s_\s(v^j) =\left(\|v^j_k\|_q \frac{v^j_{k\s(k)}}{|v^j_{k\s(k)}|}\right)_{k\le n} \in \ell_2$$
(with the convention $0/0 =0$) are orthonormal.
\end{itemize}
\end{theorem}

\begin{corollary} \lb{ge}
If  $Y \subset \ell_p(\ell_q)$, $1 < p, q < \infty$, $p \neq q$, is
$1$-complemented in $\ell_p(\ell_q)$ \wtw $Y$ is isometric to $\sum\oplus_p Y_i$, where each
$Y_i$ is isometrically isomorphic to $\ell_q^{d_i}$, $d_i = \dim Y_i \in \bbN \cup\{\infty\}$.
\end{corollary}

\begin{proof}[Proof of Corollary~\ref{ge}]
The ``only if'' part follows immediately from Theorems~\ref{lplq} and \ref{l2lq}, and the ``if'' part when $q\ne2$ is a simple consequence of \cite[Lemma~6]{Kol91} and Theorems~\ref{lplq} and \ref{l2lq}.
When $q=2$ the ``if'' part follows from Lemma~\ref{l2cor}.
\end{proof}

For the proof of Theorems~\ref{lplq} and \ref{l2lq} we use
 Theorem \ref{main} and the
following   lemmas:

\begin{lemma}
\label{disj}
Suppose that $x^1,x^2$, $\|x^1\|=\|x^2\|=1$,  are disjointly supported (in the scalar sense).
Then $\Span\{x^1,x^2\}$ is $1$-complemented in $\ell_p(\ell_q)$, with $1< p, q <\infty$, $p \neq q$, \wtw one of the following conditions holds:
\begin{itemize}
\item[(a)] 
  $\vsupp x^1 \cap  \vsupp x^2 = \varnothing $, in this case $F$ is isometric to $\ell_p^2$,
\end{itemize}

\noindent
or  
\begin{itemize}
\item[(b)] $\vsupp x^1 = \vsupp x^2$ and 
$\|x^1_i\|_q=\|x^2_i\|_q$ for all $i\in \vsupp x^1,$ in this case $F$ is isometric to $\ell_q^2$,
\end{itemize}
\end{lemma}

\begin{lemma}
\label{l2cor}
Suppose that $F \subset \ell_p(\ell_q) $ is isometric to $\ell^2_2$ and is
$1$-complemented in $\ell_p(\ell_q),$  $1\le  p, q <\infty$, $p \neq q$. Then 
\begin{itemize}
\item[(a)] $F$ is spanned by disjointly supported
vectors (in the scalar sense),
\item[(b)] 
$p=2$,
\end{itemize}

\noindent
or  
\begin{itemize}
\item[(c)] 
$q=2$ and there exists a surjective isometry
$U$ of $\ell_p(\ell_q)$ such that $UF$ is spanned by disjointly supported
vectors (in the scalar sense).
\end{itemize}
\end{lemma}

\begin{remark}
It follows immediately from Lemma~\ref{disj} that if Lemma~\ref{l2cor}$(a)$ holds then $p=2$ or $q=2$.
\end{remark}

\begin{lemma}
\label{l2iff}
Let $x, y \in \ell^n_2(\ell^m_q)$, $n,m \in \bbN\cup\{\infty\}$, $1\le   q <\infty$, $q \neq 2$, $\|x\|=\|y\|=1$.
Then $F= \Span\{x,y\}\subset\ell_2(\ell_q)$ is isometrically isomorphic to $\ell_2^2$ and $1$-complemented in $\ell_2(\ell_q)$ \wtw 
for each $i\in \vsupp x \cap
\vsupp y$ there exists a constant $C_i\in\bbC$ with
$$x_i = C_i y_i$$
and for every map $\s:\{1,\dots,n\} \lra \{1,\dots,m\}$ such that $(i,\s(i))\in\supp x \cup \supp y$ whenever $i\in \vsupp x \cup \vsupp y$,  the vectors
$$s_\s(x) =\left(\|x_i\|_q \frac{x_{i\s(i)}}{|x_{i\s(i)}|}\right)_{i\le n} \in \ell_2,   \text{ and } s_\s(y) =\left(\|y_i\|_q \frac{y_{i\s(i)}}{|y_{i\s(i)}|}\right)_{i\le n} \in \ell_2$$
(with the convention $0/0 =0$) are orthonormal.
\end{lemma}

\begin{proof}[Proof of Theorems~\ref{lplq} and \ref{l2lq}]
Since $\ell_p(\ell_q)$ is complex, if $Y$ is $1$-complemented in $\ell_p(\ell_q),$ then $Y$ has
a 1-unconditional basis $ \{v^i\}^{\dim Y}_{i=1}$.

Consider $v^i, v^j$ for some $i\ne j \le \dim Y$.

If $v^i, v^j$ belong to different Hilbert components of $Y$, then by Theorem~\ref{main}, $v^i$ and  $v^j$ are scalarly disjoint and conditions $(a), (b)$ of Theorem~\ref{lplq} follow from Lemma~\ref{disj} (since it is clear that $\Span\{v^i, v^j\}$ is 1-c
omplemented in $Y$ and therefore in $\ell_p(\ell_q)$).

If $v^i, v^j$ are in the same Hilbert component of $Y$, then $\Span\{v^i, v^j\}$ is isometric to $\ell^2_2$ and when $p\ne 2$ Lemma~\ref{l2cor}   reduces our considerations to the case of disjointly supported
vectors (in the scalar sense), where we apply Lemma~\ref{disj}.

When $p=2$ we apply Lemma~\ref{l2iff}.
\end{proof}
%page 26

\begin{proof}[Proof of Lemma \ref{disj}]
Notice that when $z\in\ell_p(\ell_q)$ then 
the norming functional $z^*$ for $z$ is given by
\begin{equation}\lb{zn}
z^* = \frac1{\|z\|^{p-1}} \sum_i \|z_i\|_q^{p-q}\sum_j |z_{ij}|^{q-1} \sgn(\overline{z_{ij}}) e_{ij}^*.
\end{equation}

By Proposition~\ref{calvert}, $\Span\{x^1,x^2\}$ is 1-complemented in $\ell_p(\ell_q)$ \wtw for all $a^1, a^2 \in \bbC$ there exist $K_s,$ ($= K_s(a^1, a^2)$), $s=1, 2$, in $\bbC$ so that
\begin{equation}\lb{cal}
(a^1x^1 +  a^2x^2)^* = K_1x^{1*} + K_2x^{2*}.
\end{equation}

This is equivalent, by disjointness of $x^1,x^2$ and \eqref{zn}, to the existence of constants $K_s,$ ($= K_s(a^1, a^2)$), $s=1, 2$, such that for all
$(i,j)\in\supp x^s$:
\begin{equation*}
\begin{split}
\frac{1}{\|a^1x^1 +  a^2x^2\|^{p-1}}  &\|a^1x^1_i +  a^2x^2_i\|_q^{p-q}  |a^sx^s_{ij}|^{q-1} \sgn(\overline{a^sx^s_{ij}}) \\ 
&= K_s(a^1, a^2) \|x^s_i\|_q^{p-q}  |x^s_{ij}|^{q-1} \sgn(\overline{x^s_{ij}}),
\end{split}
\end{equation*}
which is further equivalent to the existence of constants $C_s,$ ($= C_s(a^1, a^2)= K_s(a^1, a^2)\cdot \|a^1x^1 +  a^2x^2\|^{p-1}$), $s=1, 2$, such that for all
$i\in\vsupp x^s$:
\begin{equation*}
\frac{\|a^1x^1_i +  a^2x^2_i\|_q^{p-q}}{\|x^s_i\|_q^{p-q}} = 
\left(|a^s| + |a^t| \frac{\|x^t_i\|_q^{q}}{\|x^s_i\|_q^{q}} \right)^{\frac{p-q}{q}}
 = C_s(a^1, a^2),
\end{equation*}
where $t\ne s$ and $\{t,s\} =\{1,2\}$.

Since $p\ne q$ this is equivalent to the conditions that for all
$(i,j)\in\supp x^s$:
\begin{equation*}
\frac{\|x^t_i\|_q}{\|x^s_i\|_q} = \frac{\|x^t_j\|_q}{\|x^s_j\|_q}.
\end{equation*}

This means that either $\vsupp x^1 \cap  \vsupp x^2 = \varnothing $ or   $\vsupp x^1 = \vsupp x^2$ and 
$\|x^1_i\|_q=\|x^2_i\|_q$ for all $i\in \vsupp x^1$ (since $\|x^1\|=\|x^2\|$).

The fact that $\Span\{x^1,x^2\}$ is isometric to $\ell_p^2$ or $\ell_q^2$, resp., follows immediately.
\end{proof}

\begin{proof}[Proof of Lemma \ref{l2cor}]
\mbox{   }

First we prove that if Part $(a)$ does not hold then $p=2$ or $q=2$.
The proof follows the standard technique of showing that $\ell_2$ can be
isometrically  embedded in $\ell_p$ only when $p$ is an even integer (see
\cite{LV}) and was suggested to me by Alex Koldobsky.

Assume that $F= \Span \{x,y\}$, where $\|x\| = \|y\| = 1$ and for all
$a, b \in \C$
$$
\|ax + by\| = (|a|^2 + |b|^2)^{\frac{1}{2}}.
$$
Then
\begin{equation}
\label{diff}
(|a|^2 + |b|^2 )^{\frac{p}{2}} = \|ax + by\|^p
= \sum^n_{i=1} \biggl(\sum^m_{j=1} |ax_{ij} + by_{ij}|^q\biggr)^{\frac{p}{q}}.
\end{equation}
For each $i \leq n$  let $S_i = \{ j\le m \ : \ (x_{ij}, y_{ij})\ne (0,0)\}= \{ j\le m \ : (i,j)\in \supp x\cup \supp y\}$.
We define an equivalence relation $R_i$ on the
set $S_i$ by the condition that
%page 33
$(j_1, j_2) \in R_i $ if and only if the pairs $(x_{ij_1}, y_{ij_1})$ and
$(x_{ij_2}, y_{ij_2})$ are proportional.

Let $\Lambda_i$ denote the set of equivalence classes for $R_i$ and let
$J= \{i\le n \ : \ \card(\La_i) = 1\}$ and $M = \{i\le n \ : \ \card(\La_i) > 1\}$. For each $i\in J$ let $j_i$ be the  representative of the
equivalence class of $R_i$ (i.e. $j_i \in S_i$) and for $i\in M$ let
$\{j_\alpha\}_{\alpha \in \Lambda_i}$ be the set of representatives of each
equivalence class. Then there exist positive  constants $\{K_i\}_{i\in J}$, $\{C_{i,\alpha}\}_{i\in M, \al\in \La_i}$ so that
\eqref{diff} can be written as
\begin{equation}
\label{diff1}
(|a|^2 + |b|^2)^{\frac{p}{2}} = \sum_{i\in J} K_i |a x_{ij_i } + b y_{ij_i } |^p+
\sum_{i\in M} \biggl(\sum_{\alpha \in \Lambda_i}
C_{i,\alpha} |a x_{ij_\alpha } + b y_{ij_\alpha } |^q\biggr)^{\frac{p}{q}},
\end{equation}
where the pairs $(x_{ij_\alpha}, y_{ij_\alpha})_{\alpha \in \Lambda_i}$
are mutually linearly independent for each $i\in M$. If there exists $i\in M$ and $ \beta \in
\Lambda_i$ with $x_{ij_\beta} \cdot y_{ij_\beta} \neq 0$, then there
exist $a_0, b_0 \in \C$ with
\begin{align*}
a_0 x_{ij_\beta} + b_0 y_{ij_\beta} &= 0\\\noalign{\vskip-6pt}
\intertext{and}\noalign{\vskip-6pt}
a_0 x_{ij_\alpha} + b_0 y_{ij_\alpha} &\neq 0 \quad \text{ for all }
\alpha \neq \beta, \al\in\La_i.
\end{align*}
%page 34
Fix $b = b_0$ and differentiate \eqref{diff1} with respect to $a$ along the real
axis at $a_0$.  The left-hand side of \eqref{diff1} can be
differentiated (in this fashion) infinitely many times.  But, if $q$ is
not an even integer then the $([q] + 1)$-st derivative of
$|ax_{ij_\beta}+b_0 y_{ij_\beta}|^q$ does not exist at $a_0$.

Hence, if $q$ is not even, $x_{ij} \cdot y_{ij} = 0$ for
all $i\in M,\  j\in S_i$. In particular, this implies that $\card(\La_i) = 2$ for all $i\in M$.
% i.e.,  $x$ and $y$ are (scalarly) disjoint.

Similarly, if there exists $i\in J$  with $x_{ij_i} \cdot y_{ij_i} \neq 0$, then there
exist $a_0, b_0 \in \C$ with
$
a_0 x_{ij_\beta} + b_0 y_{ij_\beta} = 0
$, we fix $b = b_0$ and differentiate \eqref{diff1} with respect to $a$ along the real
axis at $a_0$. As before, the left-hand side of \eqref{diff1} can be
differentiated (in this fashion) infinitely many times.  But, if $p$ is
not an even integer then the $([p] + 1)$-st derivative of
$|ax_{ij_i}+b_0 y_{ij_i}|^p$ does not exist at $a_0$.

Hence, if $p$ is not even, $x_{ij} \cdot y_{ij} = 0$ for
all $i\in J,\  j\in S_i$.

Further, by Lemma~\ref{l2}, if $x^* , y^*$ denote norming functionals
for $x, y$, respectively, then $\Span \{x^* , y^*\}$ is isometric to
$\ell_2^2$ and is 1-complemented in $(\ell_p (\ell_q ))^* =
\ell_{p'}(\ell_{q'})$, where $1/p + 1/p' = 1=1/q+ 1/q'$. By \eqref{zn} relations analogous to $R_i$ defined using coefficients of $x^*$ and $y^*$ are identical as the original relations $R_i$. Hence, if $q'$
is not even, then $x^*_{ij} \cdot y^*_{ij} = 0$ for
all $i\in M,\  j\in S_i$ and if $p'$
is not even, then $x^*_{ij} \cdot y^*_{ij} = 0$ for
all $i\in J,\  j\in S_i$.

Thus, by \eqref{zn}, if $q$, $q'$, $p$, $p'$ are not all even integers, that is, $q \ne 2$ and $p \ne 2$ then $x$
and $y$ are (scalarly) disjoint.
\end{proof}

We postpone the proof of Part~$(c)$, because the proof of Lemma~\ref{l2iff} is the direct continuation of just presented argument and uses the same notation.

\begin{proof}[Proof of Lemma~\ref{l2iff}]
\mbox{  }

The ``only if'' part:
By the first part of the proof of Lemma~\ref{l2cor}, since $q \ne 2$, $x_{ij} \cdot y_{ij} = 0$ for
all $i\in M,\  j\in S_i$. We will show that $M = \varnothing$.

Since $\Span\{x,y\} = \ell_2^2$ we get, by Lemma~\ref{l2},  that for all   $a, b \in \C$ with $|a|^2 + |b|^2 = 1$, the functional $\bar{a} x^* +
\bar{b} y^*$ is norming for $ax + by$.
That is, by \eqref{zn},
\begin{equation*}
\begin{split}
 \sum_i \|ax_i&+by_i\|_q^{2-q}\sum_j |ax_{ij}+by_{ij}|^{q-1} \sgn(\overline{ax_{ij}+by_{ij}}) e_{ij}^* = \\
&\bar a \sum_i \|x_i\|_q^{2-q}\sum_j |x_{ij}|^{q-1} \sgn(\overline{x_{ij}}) e_{ij}^* + \bar b \sum_i \|y_i\|_q^{2-q}\sum_j |y_{ij}|^{q-1} \sgn(\overline{y_{ij}}) e_{ij}^*.
\end{split}
\end{equation*}

Thus if $i\in M$, by disjointness of $x_i, y_i$, for each $j$ with $(i,j)\in \supp x$ we have
\begin{equation*}
  \|ax_i+by_i\|_q^{2-q} |ax_{ij}|^{q-1} \sgn(\overline{ax_{ij}}) = 
\bar a  \|x_i\|_q^{2-q} |x_{ij}|^{q-1} \sgn(\overline{x_{ij}}).
\end{equation*}

Hence, since $q\ne 2$, 
$$  \|ax_i+by_i\|_q = 
| a | \|x_i\|_q = \|ax_i\|_q. $$
Since $x_i$ and $y_i$ are disjoint and $\ell_q$ is strictly monotone we conclude that $y_i = 0$. But this implies that $\card(\La_i) = 1$, which contradicts the fact that $i\in M$. Thus $M= \varnothing$.

Hence $J= \vsupp x\cup \vsupp y$ and
for each $i\in \vsupp x \cap
\vsupp y$ there exists a constant $C_i\in\bbC$ with
$$x_i = C_i y_i.$$

Now \eqref{diff}, \eqref{diff1} and definition of $J$ imply that  for every map $\s:\{1,\dots,n\} \lra \{1,\dots,m\}$ such that $(i,\s(i))\in\supp x \cup \supp y$ whenever $i\in \vsupp x \cup \vsupp y$  
we have (with the convention $0/0 =0$):
\begin{equation*}
\begin{split}
|a|^2 + |b|^2  &= \|ax + by\|^2
= \sum_{i\in J} \biggl(\sum_{j\in S_i} |ax_{ij} + by_{ij}|^q\biggr)^{\frac{2}{q}}\\
&= \sum_{i\in J} \left(\sum_{j\in S_i} \left(\frac{|x_{ij}|}{|x_{i\s(i)}|}|ax_{i\s(i)} + by_{i\s(i)}|\right)^q\right)^{\frac{2}{q}}\\
&= \sum_{i\in J} \left(\sum_{j\in S_i} \frac{|x_{ij}|^q}{|x_{i\s(i)}|^q}\right)^{\frac{2}{q}}|ax_{i\s(i)} + by_{i\s(i)}|^2\\
&= \sum_{i\in J}  \frac{\|x_{i}\|_q^2}{|x_{i\s(i)}|^2} |ax_{i\s(i)} + by_{i\s(i)}|^2\\
&= \sum_{i\in J}   \left|a\frac{\|x_{i}\|_q}{|x_{i\s(i)}|}x_{i\s(i)} + b\frac{\|y_{i}\|_q}{|y_{i\s(i)}|}y_{i\s(i)}\right|^2.
\end{split}
\end{equation*}
Thus the vectors
$$s_\s(x) =\left(\|x_i\|_q \frac{x_{i\s(i)}}{|x_{i\s(i)}|}\right)_{i\le n} \in \ell_2,   \text{ and } s_\s(y) =\left(\|y_i\|_q \frac{y_{i\s(i)}}{|y_{i\s(i)}|}\right)_{i\le n} \in \ell_2$$
(with the convention $0/0 =0$) are orthonormal.

The ``if'' part:
It is clear from above calculations that if $x, y$ are of the form described in Lemma~\ref{l2iff} then $\Span\{x, y\}$ is isometrically isomorphic to $\ell_2^2.$

Further, by \eqref{zn} we have for all $a, b\in\bbC$: 
\begin{equation*}
\begin{split}
 (ax &+by)^* = \frac1{(|a|^2 + |b|^2 )^{\frac12}} \sum_i \|ax_i+by_i\|_q^{2-q}\sum_j |ax_{ij}+by_{ij}|^{q-1} \sgn(\overline{ax_{ij}+by_{ij}}) e_{ij}^* \\
&=\frac1{(|a|^2 + |b|^2 )^{\frac12}} \biggl(\sum_{i\in\vsupp x\setminus\vsupp y} \|ax_i\|_q^{2-q}\sum_j |ax_{ij}|^{q-1} \sgn(\overline{ax_{ij}}) e_{ij}^* \\
&\hspace{5mm} + \sum_{i\in\vsupp x\cap\vsupp y} \|ax_i+bC_ix_i\|_q^{2-q}\sum_j |ax_{ij}+bC_ix_{ij}|^{q-1} \sgn(\overline{ax_{ij}+bC_ix_{ij}}) e_{ij}^* \\
&\hspace{5mm} + \sum_{i\in\vsupp y\setminus\vsupp x} \|by_i\|_q^{2-q}\sum_j |by_{ij}|^{q-1} \sgn(\overline{by_{ij}}) e_{ij}^* \biggr)
\end{split}
\end{equation*}
\begin{equation*}
\begin{split}
 &=\frac1{(|a|^2 + |b|^2 )^{\frac12}} \biggl(\bar a \sum_{i\in\vsupp x\setminus\vsupp y} \|x_i\|_q^{2-q}\sum_j |x_{ij}|^{q-1} \sgn(\overline{x_{ij}}) e_{ij}^* \\
&\hspace{5mm} + \overline{a+bC_i}\sum_{i\in\vsupp x\cap\vsupp y} \|x_i\|_q^{2-q}\sum_j |x_{ij}|^{q-1} \sgn(\overline{x_{ij}}) e_{ij}^* \\
&\hspace{5mm} + \bar b  \sum_{i\in\vsupp y\setminus\vsupp x} \|y_i\|_q^{2-q}\sum_j |y_{ij}|^{q-1} \sgn(\overline{y_{ij}}) e_{ij}^* \biggr)\\
&=\frac1{(|a|^2 + |b|^2 )^{\frac12}} ( \bar a x^* + \bar b y^*).
\end{split}
\end{equation*}

Thus, by Lemma~\ref{calvert}, $\Span\{x, y\}$ is 1-complemented.
\end{proof}

\begin{proof}[Proof of Lemma~\ref{l2cor}$(c)$]
Assume that $q=2$.
First notice that, if $z \in \ell_p (\ell_2 )$ and $\|z\| = 1$, then   
the norming functional $z^*$ for $z$ is given by
\begin{equation}
\label{zstar}
z^*=\sum_{i=1}^n\|z_i\|_2^{p-2} \cdot \sum _{j=1}^m
\overline{z_{ij}}e_{ij}
\end{equation}

Let $F  = \Span   \{x, y\}$, where $\|x\| = \|y\| = 1$ and
$\|ax + by\| = (|a|^2 + |b|^2 )^{1/2}$
for all $a, b \in \C$.
Let $U_i : \ell_m^2  \to \ell_m^2$  be a surjective isometry such that, for
all $i \le n$,
$$
U_ix_i = \|x_i\|_2 e_1 .
$$
Define an isometry $U  : \ell_p (\ell_2 ) \to \ell_p (\ell_2 )$ by
$$
U ((z_i)^n_{i=1}) = (U_iz_i)^n_{i=1}.
$$
Then $U F = \Span\{U x, U y\}$ is isometric to $\ell_2^2$ and is 1-complemented
in $\ell_p (\ell_2 )$.

Thus, by Lemma~\ref{l2}, for all $a, b \in \C$ with $|a|^2 + |b|^2 =
1$, we have
$$
(a U x + bU y)^*   = \bar a (U x)^*  + \bar b(U y)^*.
$$
Hence, by \eqref{zstar} and by the choice of $U$, we get for all $i
\le n$ and all $j\geq 2$:
\begin{equation}
\label{isol21}
\begin{split}
   \|a(U x)_i + b(U y)_i\|^{p-2}_2  &\cdot  \left(\bar a  \cdot 0 + \bar
b \cdot \overline{ (U y)_ {ij}}\right) = \\ 
&\bar a \cdot
  \|(U x)_i\|_2^{p-2} \cdot    0 + \bar b \cdot  \|(Uy)_i\|_2^{p-2}
\overline{(U y)_{ij}}.
\end{split}
\end{equation}

Now for each $i$ with $\|(Ux)_i\|_2\ne 0$ we consider two cases: either
\begin{enumerate}
\item there exists $j\geq  2$ with $(U y)_
{ij}\ne0$,  or
\item $(U y)_i = ((U y)_{i1}, 0, 0,\dots,0)$.
\end{enumerate}

In case $(1)$ we have
$$
\|a(U x)_i + b(U y)_i\|_2^{p-2}   = \|(U y)i\|_2^{p-2}
$$
whenever $|a|^2 + |b|^2 = 1$.  In particular, since $p\ne 2$ and $\|(U
x)_i\|_2 =\|(U y)_i\|_2$, we get
$$
\left\| a \frac{(Ux)_i}{\|(Ux)_i\|_2}+ b
\frac{(Uy)_i}{\|(Uy)_i\|_2}\right\|_2=1.
$$
Thus
${(Ux)_i}/{\|(Ux)_i\|_2}$ and
${(Uy)_i}/{\|(Uy)_i\|_2}$ form
an orthonormal basis for $\ell_2$.  Since  $(Ux)_i/\|(Ux)_i\|= (1, 0,
0,\dots)$, the vectors $(U x)_i$ and $(U y)_i$
are disjoint.

In case $(2)$, i.e. when $(U y)_i = ((U y)_{i1}, 0, 0,\dots,0)$, then by \eqref{zstar} and by
the form of $U$ we get for all $a, b$ with $|a|^2 + |b|^2 = 1$:
\begin{equation}
\label{isol22}
\begin{split}
   |a(U x)_{i1} + b(U y)_{i1}|^{p-2}&\cdot
     (\bar a (U x)_{i1} + \bar b\,\overline{(U y  )_{i1}}) = \\
&\bar a (U x)^{p-2}_{i1} (U x)_{i1} + \bar b|(U y)_{i1}|^{p-2} \overline{  (U y)_{i1}}.
\end{split}
\end{equation}

Let $a \in\R_+$ and $b = c e^{i\theta}$, where $c \in\R_+$, $a^2 + c^2  = 1$ and $\theta$
is such that $e^{-i\theta}\cdot \overline{ (U y)_{i1}} = |(U y)_{i1}|$.
Set $\alpha = (U x)_{i1} > 0$, $\beta = |(U y)_{i1}|\geq 0$.  Then
\eqref{isol22} becomes
\begin{equation}
\label{isol23}
(a\alpha + c\beta)^{p-1}    = a\alpha^{p-1}    + c\beta^{p-1}   ,
\end{equation}
and this equation holds for all $a, c \geq 0$ with $a^2 + c^2  = 1$.

For any fixed $u, v \geq 0$ define a function $f_{u,v} : [0, 1] \to \R$
by
$$
f_{u,v}(a) = au + \sqrt{ 1- a^2 } v.
$$
It is not difficult to check that $f_{u,v}$ attains its maximum on
$[0,1]$ at the point $a_0 = a_0 (u, v) = u (u^2 + v^2 )^{-1/2}$ and the
maximum value of $f_{u,v}$ is equal to $M (u, v) = (u^2 + v^2 )^{1/2}$.

Since equation \eqref{diff1}   can be written as
$$
(f_{\alpha,\beta}(a))^{p-1}=f_{\alpha^{p-1},\beta^{p-1}(a)},
$$
we have
\begin{align*}
a_0(\alpha,\beta)&=a_0(\alpha^{p-1},\beta^{p-1}),\\
M(\alpha,\beta)^{p-1}&=M(\alpha^{p-1},\beta^{p-1}).
\end{align*}
Thus
$$
\biggl(\frac{\alpha}{\sqrt{\alpha^2+\beta^2}}\biggr)^{p-1}=\frac
{\alpha}{\sqrt{\alpha^2+\beta^2}}.
$$

Since $p\neq2$, we conclude that either $\alpha = 0$ or $\alpha
(\alpha^2 + \beta^2 )^{-1/2} = 1$, i.e., $\beta = 0$.  But $\alpha =
(U x)_{i1}= \|(Ux)_i\|_2$ and $\beta = |(U y)_{i1}|= \|(Uy)_i\|_2$.  Thus $(U x)_i = 0$ or $(U y)_i =
0$.  Hence $U x$ and $U y$ are disjoint.
\end{proof}

\begin{remark}
Lemmas \ref{disj}-\ref{l2iff} are all valid (with the presented proofs) both
in the complex and real case.  We suspect that Theorems~\ref{lplq} and \ref{l2lq}, too, are true
in the real case, but our method of proof does not work then.
\end{remark}

%page 42

\section{$1$-complemented disjointly supported subspaces of Orlicz and
Lorentz spaces}
\label{gen}
In this section we fully characterize  subspaces of
(real or complex) Orlicz and Lorentz spaces that are spanned by
disjointly supported elements and $1$-complemented.

In particular, it follows from Theorems \ref{orlicz} and
\ref{glorentz} that in ``most'' Orlicz and Lorentz spaces the only
$1$-complemented disjointly supported subspaces are those spanned by a
block basis with constant coefficients (of some permutation of the
original basis). %
%page 43
\begin{theorem}
\label{orlicz}
Let $\ell_\phi$ be a (real or complex) Orlicz space and let $x,y \in
\ell_\phi$, be disjoint elements such that $\|x\|_\phi = \|y\|_\phi =
1$ and $\Span \{x,y \} $ is $1$-complemented in $\ell_\phi$.

Then one of three possibilities holds:
\begin{itemize}
\item[(1)] $\card (\supp x) < \infty$ and $|x_i| = |x_j|$
for all $i,j \in \supp x$; or
\item[(2)] there exists $p$, $1 \leq p \leq \infty$, such that
$\phi(t) = Ct^p$ for all $t \leq \|x\|_\infty$; or
\item[(3)] there exists $p$, $1 \leq p \leq \infty$, and constants
$C_1, C_2, \gamma \geq 0$ such that
$C_2 t^p \leq\nobreak \phi(t) \leq C_1 t^p$
for all $t \leq \|x\|_\infty$
and  such that, for all $j \in \supp x$,
$$|x_j| = \gamma^{k(j)} \cdot \|x\|_\infty$$ for some
$k(j) \in \Z$.
%page 44
\end{itemize}
\end{theorem}
%page 45
For the proof of the theorem we will need the following (well-known?)
lemma, whose proof is outlined in \cite{Z76}. For the convenience of
the reader we provide its proof below. %
%page 46
\begin{lemma}
\label{z}
Let $a > 0$ and suppose that $\phi: [0,a] \to \R $ is an increasing
differentiable function with $\phi(0)=0$. Suppose that there exist $a
< \alpha, \beta < 1$ so that, for all $u \leq a$,
\begin{equation}
\label{mult}
\phi(\alpha u) = \beta \phi(u).
\end{equation}
Then there exist $p > 0$ and $C_1, C_2 > 0$ such that, for all $u \leq
a$,
$$
C_2 u^p \leq \phi(u) \leq C_1 u^p.
$$
Moreover, if $\phi(u) \not\equiv C \cdot u^p$, there exists $\gamma > 0$
such that \eqref{mult} is satisfied (with the corresponding $\beta$) if
and only if $\alpha = \gamma^k$ for some $k \in\Z$.
\end{lemma}
%page 47
\begin{proof}[Proof of Theorem \ref{orlicz}]
Let $z \in \ell_\phi$. By \cite{GrzH} the norming functional $z^*$ of $z$
is given by:
$$
z^*_i = \frac1C \sgn (\overline{z_i}) \phi'
\Bigl(\frac{|z_i|}{\|z \|_\phi}\Bigr),
$$
where $C$ is a constant depending on $z$.
By Proposition \ref{conj}, for all $b \in \C$ there exist constants
$K_1, K_2$ such that
$$
(x + by)^* = K_1 x^* + K_2 y^*.
$$
Since $x$ and $y$ are disjoint, there exists a constant $K=K(b)$ so
that for all $i \in \supp x$
$$
\sgn (\overline{x_i}) \phi' \left(\frac{|x_i|}{\| x + by
\|_\phi}\right )
= K \cdot \sgn (\overline{x_i})\phi' (|x_i|).
$$
%page 48
%4 hours total
Now for all $, 0<t<1$ there exists $b \in\C$ so that $\|x + by \|_\phi
= t^{-1}$.  Thus, for all $0<t \leq 1$, there exists $C_t > 0$ so that
for all $i \in \supp x$
$$
\phi'(|x_i|\cdot t) = C_t \phi'(|x_i|).
$$
Hence for all $i,j \in \supp x$ and for all $t \leq 1$
$$
\frac{\phi'(|x_i|\cdot t)}{\phi'(|x_i|)} =
\frac{\phi' (|x_j|\cdot t)}{\phi' (|x_j|)}.
$$
Set
$$
\beta = \frac{ \phi' (|x_i|)}{\phi' (|x_j|)},
\quad u= |x_j|\cdot t ,\quad \alpha = \frac{|x_i|}{|x_j|}.
$$
In this notation we have 
$$\phi' (\alpha u)=\beta
\phi'(u)$$
 for all $u$ such that $0 \leq u \leq |x_j|$. Thus
$$\phi(\alpha u) = \beta \alpha \phi(u)$$
 for all $u$ such that $0 \leq u \leq
|x_j|$. %
%page 49

Let $j \in \supp x$ be such that $|x_j| = \|x\|_\infty$. If there
exists $i \in \supp x $ with $|x_i| \neq |x_j| = \|x\|_\infty$ then,
by Lemma \ref{z}, condition (2) or (3) holds.
\end{proof}

%page 50
\begin{proof}[Proof of Lemma \ref{z}]
Let $p= \log_\alpha (\beta)$. Let $m_0 \in \Z$ be the smallest integer
with $\alpha^{m_0} \leq a$.

If $\alpha^{m_0} < u \leq a$ we have
$$
\phi(u) \leq \phi(a) \leq \frac{\phi(a)}{\alpha^{m_0 p}} \cdot u^p
$$
and
$$
\phi(u) \geq \phi(\alpha^{m_0}) \geq
\frac{\phi(\alpha^{m_0})}{a^p}\cdot u^p.
$$
If $\alpha^{m+1} < u \leq \alpha^m$ for some $m \leq m_0$ we have
$$
\thickmuskip=2mu plus 2mu
\phi(u) \leq \phi(\alpha^m) = \beta\phi(\alpha^{m-1}) = \cdots
=\beta^{m-m_0}\phi(\alpha^{m_0})
=\frac{\phi(\alpha^{m_0})}{\beta^{m_0+1}}(\alpha^p)^{m+1} \leq
\frac{\phi(\alpha^{m_0})}{\beta^{m_0+1}} \cdot u^p
$$
%page 51
and
$$
\thickmuskip=1mu plus 2mu
\phi(u) \geq \phi(\alpha^{m+1}) = \beta \phi(\alpha^m) = \dots
=\beta^{m+1-m_0} \phi(\alpha^{m_0})
=\frac{\phi(\alpha^{m_0})}{\beta^{m_0-1}}\cdot(\phi^p)^m \geq
\frac{\phi(\alpha^{m_0})}{\beta^{m_0-1}}u^p.
$$
Set
$
C_1 = \max \{ {\phi(a)}/{\alpha^{m_0 p}}, \
{\phi(\alpha^{m_0})}/{\beta^{m+1}}\}
$
and
$
C_2 = \min \{ {\phi(\alpha^{m_0})}/{a^p},\
{\phi(\alpha^{m_0})}/{\beta^{m_0-1}}\}.
$
Then $$C_2 u^p \leq \phi(u) \leq C_1 u^p $$ for all $u$ with $0 \leq u \leq a$.

Further define a function $h_\phi: (-\infty, \ln a] \to \R$ by
$$
h_\phi(t) = \frac{d}{ds}(\ln(\phi(e^s)))\big|_{s=t}.
$$
Then, by \eqref{mult}, $h_\phi(t+ \ln(\alpha)) = h_\phi(t)$ for all $ t
\leq \ln a$.
%page 52
Thus, since $\alpha \neq 1$, either
\begin{itemize}
\item  $h_\phi$ is constant, that is, there exists a
constant $K$ so that $\phi(u) = K \cdot u^p$ for all $ u \leq a$, or
\item  $h_\phi$ is periodic, that is, there exists $w$, with minimal $|w|$,
so that  $h_\phi(t+w) = h_\phi(t)$ for all $t \leq \ln a$.
\end{itemize}
Thus there exists $\gamma>0$ (namely $\gamma=e^w$) and $k \in \Z$ such
that $\alpha = \gamma^k$.
\end{proof}

Our next theorem describes disjointly supported $1$-complemented
subspaces of (real or complex) Lorentz sequence spaces.
%page 53
\begin{theorem}
\label{glorentz}
Let $\ell_{w,p}$, with $1 < p < \infty$, be a real or complex Lorentz
sequence space.  Suppose that $\{ x_i \}_{i \in I}$ are mutually
disjoint elements of $\ell_{w,p}$ such that $\card (I) \geq 2$ and $F
= \overline{\Span} \{ x_i \}_{i \in I} $ is $1$-complemented in
$\ell_{w,p}$.  Suppose, moreover, that $w_\nu \neq 0$ for all $\nu
\leq \Sigma \defeq \sum_{i \in I} \card (\supp x_i)\quad (\leq \infty)$.

Then 
\begin{itemize}
\item[$(a)$] $w_\nu = 1$ for all $ \nu \leq \Sigma$, 
\end{itemize}

\noindent
or 
\begin{itemize}
\item[$(b)$] $ |x_{il}| = |x_{ik}|$ for
all $i \in I$ and all $k, l \in \supp x_i$.
\end{itemize}
\end{theorem}
%page 54
\begin{proof}
With each element $z \in \ell_{w,p}$ we associate a decreasing sequence of
positive numbers $(\tilde{z}_i)^{l(z)}_{i=1}$ and ``level sets'' $A_i(z)$
defined inductively as follows:
$$
\def\eqalign#1{\null\,\vcenter{\openup\jot
  \ialign{\strut\hfil$\displaystyle{##}$&$\displaystyle{{}##}$\hfil
      \crcr#1\crcr}}\,}
\eqalign{
\tilde{z}_1 &= \| z \|_\infty, \cr
\tilde{z_2} &= \max \{ |z_j|: j \in \N \setminus A_1(z) \},}
\qquad \eqalign{A_1(z) &= \{j \in \N: |z_j| = \tilde{z}_1 \},\cr
A_2(z) &= \{ j \in \N: |z_j| = \tilde{z}_2 \},}
$$
and so on.  Note that
$l(z)$ is the largest integer such that $\tilde{z}_{l(z)} > 0$ and $\supp
z = \bigcup^{l(z)}_{i=1} A_i(z)$.
%page 55

For $i \leq l(z)$ introduce also 
$$s_0(z) = 0,\ \  s_i(z) = \sum^i_{j=1}
\card (A_i(z)),$$
 $$L_i(z) = \{ s_{i-1}(z) + 1, \dots, s_i(z) \}
\subset \N,$$
 and let $\delta_i: A_i(z) \to L_i(z)$ be a bijection.

Finally, for any set $A \subset \N $ denote by $ \cal{P}(A)$ the set
of all permutations of $A $.

In this notation we can easily describe norming functionals $z^N$ for
$z$. Namely, for each $ j $ with $ 1 \leq j \leq l(z) $, there exists a family
of coefficients $ \{ \lambda_\sigma \}_{\sigma \in \cal{P}(L_j(z))}$
such that %
%page 56
$ \lambda_\sigma \geq 0, \quad \sum_{\sigma \in \cal{P}(L_j)}
\lambda_\sigma = 1 $ and
$$
\bigl(|z^N_k|\bigr)_{k \in A_j(z)} = \left (\frac{\tilde{z}_j}{\| z \| }\right
)^{p-1} \! \sum_{\sigma \in \cal{P}(L_j(z))} \lambda_\sigma
(w_{\sigma(\delta_j (k))})_{k \in A_j(z)}.
$$
In particular, we can compute the $\ell_1$-norm of $z^N$ restricted to a level set $A_j(z)$
\begin{equation}
\label{l1n}
\left\|\bigl(|z^N_k|\bigr)_{k \in A_j(z)}\right\|_{\ell_1}
= \left (\frac{\tilde{z}_j}{\| z \| }\right
)^{p-1} \! \sum_{n\in L_j(z)} w_n.
\end{equation}

Notice that the left hand side of \eqref{l1n} does not depend on the choice of the
norming functional $z^N$ for $z$. %
%page 57

Now assume that $ l(x_1) > 1 $. We will show that $ w_\nu= 1 $ for all natural numbers
$\nu \leq \Sigma$, where $\Sigma  \defeq \sum_{i \in I} \card (\supp x_i)$ as defined above.

If
$\nu \leq \Sigma$, there exists $n \leq l_1(x_1)$ such that
$$
\nu- s_n(x_1) \leq \Sigma - \card (\supp x_1).
$$

Further, there exist $ \mu \in \N$ and $ \{ j_i \}^{\mu}_{i=2} \subset
\N $ such that $ j_i \leq l(x_i)$ for $ i=2, \dots, \mu $ and
$$\sum^{\mu}_{i=2} s_{j_i}(x_i) =: M \geq \nu - s_N(x_1).$$

Choose $ \{ a_i \}^{\mu}_{i=1} \subset \R_+$ such that,
%page 58
for all $i$ with $2 \leq i \leq \mu$,
$$
a_1 (\widetilde{{x}_{1}})_{1} > a_i (\widetilde{{x}_{i}})_{j_i} > a_1 (\widetilde{{x}_{1}})_{2}
$$
and $ a_1 (\widetilde{{x}_{1}})_{n} > a_i (\widetilde{{x}_{i}})_{j}$ for all $ j > j_i $.

To shorten the notation, set $ x = \sum^{\mu}_{i=1} a_1 x_1 $. Then
\begin{equation}
\label{lov1}
a_1 (\widetilde{{x}_{1}})_{1} = \tilde{x}_1\quad\hbox{and} \quad A_1(x) = A_1 (x_1).
\end{equation}
Moreover, there exists $ k \in \N$ with $2 < k \leq 1 + \sum^{\mu}_{i=2}
j_i$ such that, for all $ \alpha $ satisfying $ 2 \leq \alpha \leq n $,
\begin{equation}
\label{lov3}
a_1 (\widetilde{{x}_{1}})_{ \alpha } = \tilde{x}_{k +(\alpha -2)},\quad A_{k+(\alpha
-2)}(x) = A_{\alpha}(x_1),
\end{equation}
and 
$$ s_{k-1}(x) = s_1 (x_1) +M. $$ 

Thus
\begin{equation}
\label{lov4}
L_1(x) = L_1(x_1) = \{ 1, \dots, s_1(x_1) \} ,
\end{equation}
%page 59
and, for all $ \alpha $ with $ 2 \leq \alpha \leq n $,
\begin{equation}
\label{lov5}
L_{k + (\alpha -2)}(x) = M + L_\alpha (x_1) = \{ M + s_{\alpha -1}(x_1) + 1,\,
\dots,\, M+ s_{\alpha} (x_1) \} .
\end{equation}

By Proposition \ref{conj} there exist norming functionals $ x^N$ for
$x$ and $x_i^N$ for $x_i$, 
%where $i=1, \dots, \mu $, 
and constants $
K_i$, where $i=1, \dots, \mu $, such that
$$
x^N = \sum^{\mu}_{i=1} K_i x_i^N.
$$
Thus, by \eqref{l1n} and \eqref{lov1},
$$
\left ( \frac{{\tilde{x}_1}}{\|x\|}\right )^{p-1} \sum_{j \in L_1(x)} w_j =
K_1 \cdot \left (\frac{a_1 (\widetilde{{x}_{1}})_{1}}{\| x_1 \|} \right )^{p-1}\sum_{j \in
L_1(x_1)}w_j.
$$
Hence, by \eqref{lov4},
\begin{equation}
\label{lov6}
K_1 \cdot\left( \frac{\|x\|}{\| x_1 \| }\right )^{p-1} = 1.
\end{equation}
%page 60

Moreover, by \eqref{l1n} and \eqref{lov3}, for all $\alpha$ with $ 2 \leq
\alpha \leq n $ we get:
$$
\left ( \frac{\tilde{x}_{k + (\alpha -2)}}{\|x\| }\right )^{p-1}
\! \sum_{j \in L_{k + (\alpha-2)}(x)} w_j = K_1 \cdot \left (
\frac{a_1 (\widetilde{{x}_{1}})_\alpha}{\|x_1 \| }\right )^{p-1} \!
\sum_{j \in L_\alpha (x_1)} w_j.
$$
Hence, by \eqref{lov6} and \eqref{lov5},
$$
\sum^{s_\alpha(x_1)}_{j=s_{\alpha-1}(x_1) + 1} w_{j+M} =
\sum^{s_\alpha(x_1)}_{j= s_{\alpha-1}(x_1)+1} w_j.
$$
Since $\{w_j\}$ is a decreasing sequence of numbers we immediately
conclude that, for all $ \alpha$ with $2 \leq \alpha \leq n $,
$$
w_{s_{\alpha-1}(x_1)+1} = w_{s_\alpha (x_1) + M}.
$$
%page 61
Since $M \geq 1 $ we get
\begin{equation}
\label{lov7}
w_{s_1(x_1)+1} = w_{s_2(x_1) +M} = w_{s_2(x_1)+1}=w_{s_3(x_1)+M} =
\dots = w_{s_n(x_1)+M}.
\end{equation}

Finally, choose $ \{ b_i \}^{\mu}_{i=1} \subset \R_+ $ in such
a way that, for all $ i$ with $2 \leq i \leq \mu $,
$$
b_i (\widetilde{{x}_{i}})_{j_i} > b_1 (\widetilde{{x}_{1}})_{1}.
$$
%page 62
Now set $ y = \sum^{\mu}_{i=1} b_i x_i $. Then there exists $ t \in \N
$, with
$1 \leq t \leq 1 + \sum^{\mu}_{i=2} j_i$, such that for all $\alpha$
with $1 \leq \alpha \leq n $ we have
\begin{equation}
\label{lov8}
b_1 (\widetilde{{x}_{1}})_{ \alpha} = \tilde{y}_{t + (\alpha -1 )}\ \ ; 
A_\alpha(x_1) = A_{t+(\al -1)}(y).
\end{equation}
Similarly, as before,
$$
s_{t+(\alpha -1)}(y) = s_\alpha(x_1)+M
$$
and
\begin{equation}
\label{lov9}
L_{t+(\alpha-1)}(y) = M + L_\alpha (x_1) = \{ M + s_{\alpha-1}
(x_1)+1,\, \dots,\, M + s_\alpha(x_1) \}.
\end{equation}

Again, by Proposition \ref{conj} there exist norming functionals 
$y^N$ for $y$ and $x_i^N$ for $x_i$,
% where $i=1, \dots , \mu $, 
and
constants $ K_i'$, where $i=1, \dots, \mu $, such that
$$
y^N = \sum^\mu_{i=1} K_i' x_i^N.
$$
%page 63
Thus, by \eqref{l1n} and \eqref{lov8} we get, for all $ \alpha$ with $1
\leq \alpha \leq n $,
$$
\left ( \frac{\tilde{y}_{t+(\alpha-1)}}{\| y \| }\right )^{p-1} \cdot
\sum_{j \in L_{t+(\alpha-1)}(y)} w_j = K_1' \! \left (
\frac{(\widetilde{{x}_{1}})_{ \alpha}}{\|x\|}\right )^{p-1} \cdot \sum_{j \in
L_\alpha(x_1)} w_j.
$$
Hence, by \eqref{lov9},
\begin{equation}
\label{lov10}
\left (\frac{1}{\| y \| }\right)^{p-1} \cdot
\sum^{s_\alpha(x_1)}_{j=s_{\alpha-1}(x_1) + 1} w_{j + M} = K_1'
\cdot \left (\frac{1}{\|x\|}\right )^{p-1}
\sum^{s_\alpha(x_1)}_{j=s_{\alpha-1}(x_1)+1} w_j.
\end{equation}
If $ \alpha = 2 $, by \eqref{lov7} we conclude that
$$
\left (\frac{1}{\| y \|}\right )^{p-1} = K_1' \cdot \left
(\frac{1}{\|x\|}\right )^{p-1}.
$$
%page 64
Thus, when $ \alpha = 1 $, \eqref{lov10} becomes
$$
\sum^{s_1(x_1)}_{j=1} w_{j+M} = \sum^{s_1(x_1)}_{j=1}w_j.
$$
Thus 
$$ w_1 = w_{s_1(x_1) + M} $$
 and since $ M \geq 1 $, by
\eqref{lov7} we get 
$$ w_1 = w_{s_n(x_1) + M}.$$  Since $ \nu \leq
s_n(x_1) + M $ we conclude that 
$$ 1 = w_1 = w_\nu,$$ 
which ends the
proof.
\end{proof}

  \bibliographystyle{alpha}
  \bibliography{tref}

\end{document}